\magnification=\magstep1
\input amstex
\documentstyle{amsppt}
\NoRunningHeads
\pageheight{8.5truein}
\nologo
\NoBlackBoxes
\topmatter
\title Kirillov Theory, Treves Strata, Schr\"odinger Equations and
analytic hypoellipticity of Sums of Squares
\endtitle
\author Sagun Chanillo
\endauthor
\affil Rutgers, The State University\endaffil
\address Dept. of Math. 
110 Frelinghuysen Rd. 
Piscataway  NJ  08854\endaddress
\endtopmatter
\document

\subhead \S 1. Introduction\endsubhead
\medskip

Let $N$ be a connected nilpotent Lie group with associated Lie Algebra
$\frak n$ of left-invariant vector fields.  We shall assume that
$\frak n$ is stratified in the sense of Folland [F].  That is, there
is a direct sum decomposition of finite-dimensional vector spaces
$V_i$, with
$$
\frak n = \overset{m+1}\to{\underset{j=1}\to\oplus} V_j, \left[ V_j,
V_i\right] \subseteq V_{i+j}, i+ j\leq m + 1,
$$
and $\left[V_j, V_i\right]=0$ if $i+j > m+1$.  Furthermore, we assume
that $V_1$ generates the Lie algebra $\frak n$.  Let $\dim V_1 =
p$, and let $\Cal B=\{X_1, \dots, X_p\}$ be a basis for $V_1$.  We
study the PDE (and associated operator $L$), 
$$
Lu = \sum^p_{i=1} X_i^2 u =0.\tag 1.1
$$
In particular, we wish to study the analytic hypo-ellipticity of the
operator $L$.  The $C^\infty$ hypo-ellipticity of $L$ follows by a
celebrated theorem of H\"ormander [Ho].
\medskip

An important counter example of Baouendi-Goulaouic [BG] demonstrates
that $L$ need not be analytic hypoelliptic.  However, there are
situations where (1.1) can be analytic hypoelliptic, like on the
Heisenberg group and for the related operators of the type of
Baouendi-Grushin.  This was demonstrated by Treves [Tr1] using a
microlocal approach and a parametrix construction, and by Tartakoff [T]
using energy estimates.  See [Tr2] for a precise statement and a
historical overview.  After the [BG] example, several examples have
been found by Oleinik [O], Hanges-Himonas [HH1-2], Helffer [H],
Pham-The-Lai and Robert [LR] and by Christ [C1-C2].  More recently
there have been articles by Hoshiro [Hos] and the Costins' [CC].  We
will put these results in perspective in the sequel.
\medskip

The point of view we adopt is that from [Tr2].  We will thus be
working on the co-tangent bundle of $N, T^*N$, and also we will have
reason later to work on the complexified co-tangent bundle, $\Bbb C
T^* N$.  On $T^*N$, we will focus on the symbols of the vector fields
$V_1$.  The characteristic  set $\Sigma$ of (1.1) then has a
stratification given by the vanishing of the symbols of the vector
fields in $V_1$ and by the vanishing of the symbols of the vector
fields we systematically get by taking brackets of the vector fields
in $V_1$, [Tr2].  We assume that one such stratum is not symplectic.
That is to be more precise:  The equation defining the non-symplectic
stratum is given by the vanishing of a co-vector.
\medskip

Thus we are assuming that one can find a curve $\gamma(t)$ on $N$, and
more precisely in the base projection of $\Sigma$ on $N$, such that the
tangent line bundle of $\gamma (t)$ is orthogonal for the fundamental
symplectic form to $TS$ restricted to $\gamma, TS |_\gamma$, where $S$
is a Poisson stratum of $\Sigma$ that is not symplectic.  Of course it may
happen, as in the example [BG], that $\Sigma$ itself is non-symplectic.
\medskip

Our approach is via the theory of representations due to Kirillov.  We
may naturally attach elements $\lambda\in \frak n^*$ (the dual Lie
algebra) to the vanishing co-vectors that define the non-symplectic
strata.  One then induces representations for these elements
$\lambda$.  These representations give rise to Schr\"odinger equations
on some $\Bbb R^n$, $n\geq 1$.  The principal difficulty arises when
$n\geq 2$.  The Schr\"odinger equations then have potentials which in
general go off to negative infinity and are thus not confining.
Bounds for such potentials are difficult to come by,  though we ask
for very weak bounds.  Nevertheless, we are able to obtain bounds in
an important case, thereby obtaining solutions to (1.1) which are
smooth but not analytic.  Previous works on this topic when examined
under the view point in this article shows that the majority of 
the examples considered lead to
representations that act on functions on $L^2(\Bbb R)$, that is the
so-called Gelfand-Kirillov dimension is 1 ([H] for example is not
restricted to the Kirillov dimension being 1).  Therefore instead of a
true PDE Schr\"odinger, the analysis is that of an ODE Schr\"odinger.
However, our article shows that to understand completely the
conjecture (1.1) in [Tr2], the analysis of Schr\"odinger equations in
$\Bbb R^n$ with degenerate potentials is unavoidable.  Our treatment
of the Schr\"odinger equations uses machinery from potential theory, in
particular Harnack's inequality.  Thus we rely extensively on
positivity.  This is a technical drawback and prevents us from
capturing the oscillatory solutions which will necessarily arise in
the general case.  However, we state a theorem of Treves that
indicates that the bounds we desire are valid in full generality.
\medskip

Since, we are assuming that the curves $\gamma (t)$ lie in the base,
we are also ruling out the interesting example of M\'etivier [M].  The
idea of using representations in the study of (1.1) is not new and
goes back to Rothschild-Stein [RS], and subsequently was developed by
Rockland [Ro] and Helffer-Nouriggat [HN] all in the
$C^\infty$-hypoellipticity context.  The novel point here is to
associate representation theory with the notion of strata [Tr2].
\medskip

Lastly, we do not have general theorems when the vector fields are not
free.  The notion of lifting vector fields introduced in  [RS], has to
be handled with care as we will show by examples.  In any case the
bounds for the Schr\"odinger equation allow us to easily construct new
classes of sums of squares which are composed of degenerate vector
fields and which have smooth, non-analytic solutions.
\medskip

The paper is organized as follows.  In Section 2, we show how the
induced representations are linked to the co-vectors that define the
strata of [Tr2]. These give rise to Schr\"odinger operators.  It is then
shown how solutions to these Schr\"odinger equations if they satisfy a
certain growth bound, yield non-analytic solutions of (1.1).  In
Section 3, under additional conditions we obtain polynomial growth
bounds on the solutions to the Schr\"odinger equations of Section 2.
The last section, section 4, contains several examples, which
clarify and focus on what remains to be done.
\medskip

\noindent{\bf Acknowledgement}  This project started out as a joint work with
F. Treves.  It was he who suggested to us an investigation of the link
between strata and representations.  We have benefitted greatly by his
comments and encouragement as the work progressed.  I also thank
Diptiman Sen for several pertinent comments regarding Section 3.  I
also wish to thank Andrea Malchiodi for a very helpful conversation
regarding Lemma (3.7).  Finally my hearty thanks to Shri S. Devanandji
for penetrating comments on Section 3.  Supported in part by NSF grant
DMS-9970359.

\subhead \S 2.  Induced Representations and Strata\endsubhead
\medskip

We will freely use the theory of induced representations as developed
by Kirillov [K].  An excellent account is given in the book by Corwin
and Greenleaf [CG].  We assume that 

$$
 \frak n =\overset{m+1}\to{\underset{j=1}\to\oplus} V_j, \tag 2.1
$$
and that the co-vector, the vanishing of which defines the
non-symplectic stratum occurs in $V^*_{s+1}$, the dual vector space to
$V_{s+1}$ in (2.1).  This co-vector can always be picked as a basis
element for $V^*_{s+1}$.  This particular co-vector will be denoted by
$\lambda_1$, 
$$
\lambda_1 = (0, 0, \dots, 0, \lambda_{s+1,\alpha} , 0, \dots, 0 ),
\lambda_{s+1, \alpha} \neq 0.
$$
More precisely, $\lambda_1\in \frak n^*$, with $\lambda_1|_{V_i}=0$, 
$i\neq s+1$, and we select a basis $\{X_{s+1,j}\},\,  j=1,\dots,\dim V_{s+1}$ 
for $V_{s+1}$, such that in this basis $\lambda_1(X_{s+1,j})=0, 
j\neq \alpha$, and $\lambda_1(X_{s+1,\alpha})=\lambda_{s+1,\alpha}$.

Next we pick an element in $V^*_{m+1}$.  This element is typically all
or some of the non-vanishing co-vectors in the set of equations that
describe the characteristic set $\Sigma$ for the operator $L$ defined 
in (1.1) (see remark
(2.6) and section 4).  We call this element $\lambda_2$, that is
$$
\lambda_2 = (0, 0, \dots, 0, \lambda_{m+1,\beta}), \lambda_{m+1,\beta}
\neq 0.
$$
Again more precisely, $\lambda_2\in \frak n^*$, such that, $\lambda_2|_{V_i}
=0$ if $i\neq m+1$, and we select a basis for $V_{m+1}$ which we denote by,
$\{X_{m+1,j}\},\, j=1,\dots , \dim V_{m+1}$ such that in this basis
$\lambda_2(X_{m+1,j})=0$
for $j\neq \beta$, and $\lambda_2(X_{m+1,\beta})=\lambda_{m+1,\beta}.$

To keep the notation simple we have assumed that $\lambda_2$ has a
single non-zero component.  This assumption is unnecessary as can be
seen by examining the ensuing proof.  See also remark (2.6).  Now
define,
$$
\lambda = \lambda_1 + \lambda_2\tag 2.2
$$
One of the reasons for the existence of non-analytic solutions to (1.1)
seems to be $s< m$.  This fails for the Heisenberg group, and so we
have analytic solutions.  We do not aim for maximum generality, and
thus to keep the combinatorial aspects simple we will make an
assumption.  The first assumption has been discussed a few sentences
above.

\medskip

 {\bf Assumptions:}
 
$$
\align \qquad &  s< m  \qquad  \tag 2.3\\
\qquad &  \text{For } k > 1, \, j > 1, \qquad \langle  \lambda_i, [V_k, V_j]\rangle = 0, \,   i=1,2,
\tag 2.4
\endalign
$$
where $\langle, \rangle$ is the duality bracket between tangent
vectors and co-vectors.
\medskip

All the examples in [Tr2] of non-analytic hypoellipticity that ensue
from a nilpotent group situation, are seen to satisfy both
assumptions.  In fact in all the known examples ensuing from groups
[BG], [C1,2], [HH1, 2], [H], [LR], [O], [Hos], [CC] a very strong form
of (2.4) is assumed, $[V_k, V_j] = 0, \, k, j>1$, except Example (3.5)
in [Tr2] due to N. Hanges. This example satisfies (2.4) but not the
strong condition $[V_k,V_j]=0, \, k,j>1$.
\medskip

It may occur to the reader to lower the value of $m$, by modding out
the ideal $\frak g$,

$$
\frak g = \overset{m+1}\to{\underset{j=s +3}\to\oplus} V_j,
$$
from the Lie algebra in (2.1).  One can do this at a price. Modding
out $\frak g$ leads us to a new Lie Algebra $\frak p=\frak n/\frak g$,
and its associated sums of squares operator $\tilde L$. It is possible
to construct explicit examples of non-analytic solutions to $\tilde L$
that have a smaller Gevrey order than non-analytic solutions that
can be constructed for the original operator $L$ associated to
the Lie algebra $\frak n$. Thus some care has to be exercised 
when using inductive arguments based on
modding out of ideals, if the aim is also simultaneously to
obtain solutions which have optimal Gevrey properties.
\medskip

We now wish to study the representation $\pi_\lambda$, and in
particular the derived representation, $d\pi_\lambda$, and
$d\pi_\lambda (L)$.  We have, 

\proclaim{Theorem 2.1} Under the assumption (2.3), (2.4) and for
$\lambda$ as in (2.2), 
$$
\align
d\pi_\lambda (L) =& \sum^n_{j=1}\left(\frac{\partial}{\partial t_j} +
i \left(\lambda_{s+1, \alpha } \tilde p_j (t) + \lambda_{m+1,\beta }
\tilde q_j (t)\right)\right)^2\\
 -& \sum^r_{k=1} \left(\lambda_{m+1, \beta } q_k(t) + \lambda_{s+1,
\alpha } p_k (t)\right)^2.\tag 2.5
\endalign
$$
Each $(\tilde q_j(t))^n_{j=1} , \left(q_k(t)\right)^r_{k=1}$ are
homogeneous polynomials of degree $m$, with real coefficients.
$(\tilde p_j (t)^n_{j=1}$ and $(p_k (t))^r_{k=1}$ are homogeneous
polynomials of degree $s$, with real coefficients.
\endproclaim

\remark{Remarks}  We will give an explicit algorithm for the
Gelfand-Kirillov dimension $n$, and also for the number $r$.  We will
then show that under additional structural hypotheses on the Lie
Algebra $\frak n$, the magnetic potentials $\tilde p_i (t)$ and
$\tilde q_i (t)$ vanish.  This structural hypotheses is automatically
fulfilled when $n=1$, that is the ODE situation.  Thus no magnetic
potentials arise in any of the references listed.
\endremark

\demo{Proof}  We begin by giving the algorithm for $n$ and $r$.  
\enddemo

\noindent{\bf Case 1:}  $s\neq 1, m\neq 1$. We point out that $s=0$
could possibly occur in this case. Now in this case $\lambda_1 |_{V_2}=0$.
Let $\Cal B$ be the basis for $V_1$ selected to define the operator
$L$.

Define

$$
\Cal S_1 = \left\{X_k \in \Cal B : \exists  Y_k \in \frak n , \ni
\langle \lambda_1, [X_k, Y_k]\rangle \neq 0 \right\}.
$$
Note since $\lambda_1|_{V_2}=0$, we conclude $Y_k \notin V_1$.  Once
the elements of $\Cal S_1$ have been picked and set aside we start 
afresh and define

$$
\Cal S_2 = \left\{ X_k \in \Cal B:\exists \widetilde Y_k \in \frak n, \ni
\langle\lambda_2, [X_k, \widetilde Y_k]\rangle \neq 0 \right\}.
$$
Let

$$
\Cal S = \Cal S_1 \cup \Cal S_2, \# \Cal S=n. \tag 2.6
$$
Note $\Cal S \neq \phi$, since $\frak n$ is stratified.

Define

$$
r=\dim V_1 - n= p-n. \tag 2.7
$$
When $s=1, \Cal S_1$ is defined differently, but the definition for
$\Cal S_2$ remains the same, and so does the definition of $\Cal S, n$
and $r$.  
\medskip

\noindent{\bf Case 2:}  $s=1$.  We define the pair set
$$
\Cal P_0 = \left\{ \{X_k, X_j\} \big| X_k, X_j \in \Cal B,
\ni\langle\lambda_1, [X_k, X_j]\rangle \neq 0\right\}.
$$
Define $X_{k_1}$ to be any vector field, that appears in the most
number of pairs (elements of $\Cal P_0$) in $\Cal P_0$.  There could
be several candidates, so simply pick any one.
\medskip

Then $X_{k_1} \in \Cal S_1$.  From $\Cal P_0$, eliminate all pairs
containing $X_{k_1}$.  We get a new pair set $\Cal P_1$, with 
$$
\Cal P_1 \subsetneqq \Cal P_0.
$$
We repeat the process, and from the pairs in $\Cal P_1$, let
$X_{k_2}$ be that vector field that appears in the most pairs of $\Cal
P_1$.  Then $X_{k_2}\in \Cal S_1$.  Delete all pairs in $\Cal P_1$
which contain $X_{k_2}$, we arrive at $\Cal P_2, \Cal P_2 \subsetneqq
\Cal P_1$.  Containing this way, we arrive at $\Cal P_j$, where each
pair contains vector fields that appear only once.  Pick from each
pair any one vector field and add it to $\Cal S_1$.  The process now
stops and we have $\Cal S_1$.  $\Cal S_2$ is defined as in Case 1, and
$\Cal S= \Cal S_1 \cup \Cal S_2$.
\medskip

\noindent{\bf Case 3:}  $m=1$.  Now necessarily $s=0$, and $\Cal S_1 = \phi$.
We repeat the construction of Case 2, with $\lambda_1$ replaced by
$\lambda_2$ in every instance, to construct $\Cal S_2$
\medskip

Let, $W$ be the vector subspace of $V_1$ defined by, taking $\Cal
B\backslash \Cal S$ as a basis.  Form the vector subspace $\frak h$ of
$\frak n$, 
$$ \frak h=
W\oplus
\overset{m+1}\to{\underset{j=2}\to\oplus} V_j.\tag 2.8
$$
\noindent{\bf Claim: }  \, $\frak h$ is a maximal sub-ordinate sub-algebra
for $\lambda$ given by (2.2).

It is clear from the construction that $[\frak h, \frak h] \subseteq
\frak h$, and thus $\frak h$ is a sub-algebra.  In fact $[\frak n,
\frak h] \subseteq \frak h$ so $\frak h$ is an ideal.  We now verify,
$$
\langle \lambda, [ \frak h, \frak h]\rangle = 0 \tag 2.9
$$
(2.4) will be used now to simplify the argument.  Let $h_i \in \frak
h, i =1, 2$.  Then,
$$
h_i = w_i + x_i \, , \, w_i \in W, x_i \in \underset{j\geqslant
2}\to\oplus V_j.
$$
Thus,
$$
[h_1, h_2] = [w_1, x_2] + [w_2, x_1] + [w_1, w_2] + [x_1, x_2].
$$
Now (2.4) asserts that
$$
\langle \lambda, [x_1, x_2]\rangle = 0.
$$
We have,
$$
\langle \lambda, [ h_1, h_2] \rangle = \langle \lambda, [ w_1,
w_2]\rangle + \langle \lambda, [w_2, x_1]\rangle + \langle \lambda, [
w_1, x_2]\rangle.
$$
\medskip

{\bf Case 1: } $s\neq 1,m\neq 1$.  Since $[w_1, w_2] \in V_2$, the first term on the
right vanishes.  By the construction of $\Cal S$ and $W$ the second
and third terms on the right above also vanish.  We conclude (2.9).
\medskip

{\bf Case 2: }  $s=1$.  Now $[w_2, x_1], [w_1, x_2] \in \underset{j\geqslant
3}\to\oplus V_j$.
\medskip

By the construction of $\Cal S_2$,
$$
\langle\lambda, [w_1, x_2]\rangle = \langle \lambda, [w_2, x_1]\rangle
= 0.
$$
Again by the construction of $\Cal S_2$, in the case $s=1$, and by the
construction of $W$, we see easily
$$
\langle \lambda, [w_1, w_2] \rangle = 0.
$$

\noindent{\bf Case 3:}  $m=1$.  This case is treated exactly as in Case 2.
\medskip

$\frak h$ is maximal, since inclusion of any $X_k\in \Cal S$ in $\frak
h$, clearly violates $\langle \lambda, [\frak h, \frak h]\rangle = 0$,
by the construction of $\Cal S$.  Let $H=\exp \frak h$.  Since $\frak
h$ is an ideal, $H$ is a normal sub-group of $N$.  Our
representation will act on $C^\infty (N/H)$.  It follows from (2.6) and
(2.8) that, 
$$
\dim (N/H) = \# \Cal S = n.
$$
We will identify the homogeneous space (group) $N/H$, with $\Bbb R^n$.
We consider the one-dimensional representation on the subgroup $H,
P_\lambda,$
$$
P_\lambda( h)=e^{i\langle \lambda, h\rangle},\tag 2.10
$$
where $\exp h \in H$.  We now induce a representation on $N$, using
$P_\lambda$, to get $\pi_\lambda$,
$$
\pi_\lambda = \underset{H\uparrow N}\to{\text{ ind } } P_\lambda.
$$
The number $n$ introduced above is called the Gelfand-Kirillov
dimension of $\pi_\lambda$.  Let $f\in \Cal S (\Bbb R^n)$ in the
Schwartz class.  Let $t\in \Bbb R^n$.  Assume exponential coordinates
on $N$ have been picked so that the first $n$ components correspond
exactly to the $n$-elements of $\Cal S$.  Thus for $f\in \Cal S
(N/H)$, we may write,
$$
f(t) = f (t, 0) , t \in \Bbb R^n, (0) \in \Bbb R^d,
$$
where $d=\dim N - n$. 
From now on $a_1, \dots, a_n, a_{n+1} \dots, a_{n+r}, n+r = p$ will
denote the exponential coordinates corresponding to the basis vectors
$\Cal B$ of $V_1$.  When we write $a_j, j \geqslant p + 1$, we mean
simply exponential coordinates of the basis vectors for $V_j,
j\geqslant 2$.  The last coordinates $a_j, j\geqslant p + 1$ have an
insignificant role to play.  Also all quadratic terms $a_i a_j$ no
matter what $i, j$ and cubic terms etc. will have no role to play, so
we do not take much care in explicitly writing them down.
\medskip

Let $R_a$ be the right regular representation.  We have, for $a,b\in N$,
and $h\in C^\infty(N)$
$$
R_a h (b) = h
\left(b \bullet a\right)
$$
Since the group law is always abelian in the first $p=\dim V_1$
variables, we have, 
$$
(t, 0) \bullet a = (t_1 + a_1, \dots, t_n + a_n, a_{n+1}, \dots,
a_{n+r}, \bullet, \bullet, \bullet) \tag 2.11
$$
where the remaining components $\bullet$ are to be computed using
the

\noindent 
Campbell-Hausdorff formula, and $r$ is given by (2.7).  We re-write
the right side of (2.11) as 
$$
(0, \dots, 0, a_{n+1},\dots, a_{n+r}, \bullet, \bullet, \bullet)
\bullet (t_1 + a_1, \dots, t_n + a_n, 0)\tag 2.12
$$
where the component $(0)$ in the first term in (2.11) is a vector with
$n$ components.  Again the components $\bullet$ are given by the
Campbell-Hausdorff formula.  What matters now are the $s+1,\alpha$
component and $m+1, \beta$ component in the first vector in (2.12),
since we have to use (2.10).  We have, in (2.12),
$$
\align
\text{(component)}_{s+1,\alpha} =& \sum^n_{j=1} \tilde p_j(t) a_j +
\sum^r_{j=1} p_j(t) a_{j+n}\\
 +& \sum_{j\geqslant p + 1} b_j (t) a_j + \sum p_{j_1j_2\dots
j_k} (t) a_{j_1} a_{j_2} \cdots a_{j_k}.\tag 2.13
\endalign
$$
The natural gradation on $\frak n$ shows that $\tilde p_j (t)$ and
$p_j(t)$ are homogeneous polynomials of degree $s$.  Similarly we have
in (2.12), 
$$
\align
\text{(component)}_{m+1,\beta} =& \sum^n_{j=1} \tilde q_j (t) a_j +
\sum^r_{j=1} q_j(t) a_{j+n}\\
 +& \sum_{j\geqslant p + 1} b_j' (t) a_j + \sum q_{j_1j_2\dots
j_k} (t) a_{j_1} a_{j_2} \dots a_{j_k}.\tag 2.14
\endalign
$$
Again the polynomials $\tilde q_j (t)$ and $q_j(t)$ are homogeneous of
degree $m$, with real coefficients.  Thus from (2.10), (2.12) - (2.14),
we have  for $i=\sqrt{-1}$,
$$
\align
\pi_\lambda (a) f(t) = & \exp \left( i \lambda_{s+1,
\alpha}\left(\sum^n_{k=1} \tilde p_k (t) a_k + \sum^r_{j=1} p_j (t)
a_{j+n}\right)\right.\\
 +& i \lambda_{m+1, \beta} \left(\sum^n_{k=1} \tilde q_k(t) a_k +
\sum^r_{j=1} q_j (t) a_{j+n}\right)\\
  +& h.o.t^{ }_{ } \Biggr) f(t_1 + a_1, \dots, t_n +
a_n)\tag 2.15
\endalign
$$  
where $h.o.t$ stands for the higher order terms in (2.13) and (2.14),
which eventually will not count towards the derived representation.

In fact in (2.15) we only display and will only need the linear terms
in $a_j$, for $1\leq j \leq p = n+r$.  Our goal is to compute for
$t\in \Bbb R^n$,
$$
\frac{\partial \pi_\lambda(a)}{\partial a_j} f (t) \bigg|_{a=0} , j=1, 2,
\dots, p.\tag 2.16
$$
This will allow us to compute, $
d\pi_\lambda (X_k) f(t)
$
for every $X_k\in \Cal B$, the basis for $V_1$.  From (2.15) and
(2.16) we see, for $j=1, 2, \dots, n,$
$$
\frac{\partial \pi_\lambda (a)}{\partial a_j} f (t) \bigg|_{a=0} =
\frac{\partial f}{\partial t_j} + i \left( \lambda_{s+1,\alpha } \tilde
p_j (t) + \lambda_{m+1, \beta } \tilde q_j (t) \right).\tag 2.17
$$
For $j = n+1, \cdots, n+r$, from (2.15) we have,
$$
\frac{\partial \pi_\lambda (a)}{\partial a_j} f (t) \bigg|_{a=0} =
i\left(\lambda_{s+1, \alpha } p_j (t) + \lambda_{m+1, \beta } q_j
(t)\right).\tag 2.18
$$
Thus, from (2.17), (2.18),
$$
\align
d\pi_\lambda \left(\sum^p_{k=1} X^2_k\right) f (t) =&
\sum^n_{j=1}\left(\frac{\partial}{\partial t_j} + i
\left(\lambda_{s+1, \alpha } \tilde p_j (t) + \lambda_{m+1,\beta }
\tilde q_j (t)\right)\right)^2 f(t)\\
 -&  \sum^r_{k=1} \left(\lambda_{s + 1, \alpha } p_k (t) +
\lambda_{m+ 1, \beta } q_k (t)\right)^2 f(t)
\endalign
$$ 

This establishes the theorem.

$ \hfill \blacksquare$

\proclaim{Proposition 2.2}  If the vector fields in $\Cal S$ commute
amongst themselves,
or if $\# \Cal S = 1$, then the magnetic potentials, $\tilde p_j (t),
\tilde q_j (t)$ vanish for $1\leq j \leq n$.
\endproclaim

\demo{Proof}  We again need to use the Campbell-Hausdorff formula.  We
have in (2.11), that the first $n$ components commute among
themselves.  Thus the components represented by $\bullet$
in (2.11) will have the property that those terms in $\bullet$, which
are linear in $a_j$, will only involve $a_j, n+ 1 \leq j \leq n + r$
out of the set $\{a_1, \dots, a_p\}, p= n+r$.  That is the variables
$a_1, \dots, a_n$ will not appear in the linear terms in $\bullet$.
In fact the linear terms arise on consideration of $[X, Y], X = (x,
0), Y = (a)$.  The Campbell-Hausdorff formula consists of commutators
of $[X, Y]$.  Since $[X, Z] = 0, Z= (a_1, \dots, a_n, 0)$, we see that
the mixed terms in the $j$-component of $[X, Y], j\geqslant p+1$, will
only be of the type $a_k t_j, n + 1 \leq k \leq p, 1\leq
j \leq n$.  As seen earlier it is only the terms that are linear in
$a_k, 1\leq k \leq p$, that matter in the derived representation. Thus
the $s+1, \alpha$ component in (2.11) is of the form,
$$
a_{s+1, \alpha} + \sum^r_{j=1} p_j(t) a_{j+n}  + \sum_{j\geqslant p
+1} b_j (t) a_j + h.o.t.\tag 2.19
$$
$h.o.t$ involves quadratic and higher terms in $a_j$.  We re-write
(2.19) in the form (2.12).  The $s+1,\alpha$ term in the first term in
(2.12) must certainly contain (2.19) by the Campbell-Hausdorff
formula.  To see the additional linear terms involving $a_j, 1\leq j
\leq p$, in the $s+1, \alpha$ component, we proceed by induction based
on the grading of $\frak n$.  We claim that the additional linear
terms in $a_j$ also do not involve any $a_j, 1\leq j \leq n$.  This is
certainly obvious for the terms of the lowest grade, in (2.12), the
starting point of the induction.  For the higher gradations its a
consequence of (2.19), the Campbell-Hausdorff formula and the
inductive hypothesis that any linear term in $a_j, 1\leq j \leq p$, in
a fixed vector space $V_k$ (gradation) involves only $a_j, n + 1 \leq
j \leq p$, and that there are no zero order terms in $a_j$, that is
pure polynomials in $t$.  Thus, the only additional linear terms
involving $a_j, 1 \leq j \leq p$, in the $s+1, \alpha$ component will
again be,
$$
\sum^r_{j=1} p_j (t) a_{j+n}.
$$
Thus in (2.13) $\tilde p_j (t) = 0$.  Similarly in (2.14) $\tilde q_j
(t) =0, 1 \leq j \leq n$.  Our proposition is proved.

$ \hfill  \blacksquare$
\enddemo

We now show how to construct a non-analytic solution to $L$. 

\proclaim{Theorem 2.3}  Assume for $\lambda_{s+1, \alpha} \in \Bbb C$
and $\lambda_{m+1, \beta} \in \Bbb R$, $\lambda_{m+1,\beta}\neq 0$, there exists a non-trivial,
smooth solution $f$ in all $\Bbb R^n$ to, 
$$
\align Af = & \left[\sum^n_{j=1} \left(\frac{\partial}{\partial t_j} + i
\left(\lambda_{s+1, \alpha } \tilde p_j (t) + \lambda_{m+1,\beta }
\tilde q_j (t)\right)\right)^2 \right.\\
 -& \left. \sum^r_{k=1} \left(\lambda_{s + 1, \alpha } p_k (t) +
\lambda_{m+ 1, \beta } q_k (t)\right)^2 \right] f=0.
\endalign
$$
Furthermore, assume that, there exist $C_1, C_2 > 0$, such that,
$$
|f(t) | \leq C_1 e^{C_2|t|^{s+1}}.\tag 2.20
$$
Then there exists a smooth solution $u$ to (1.1) in some neighborhood
of the origin which satisfies the following properties.
\medskip

(1)  If $f(0) = 0$, the solution $u$ to (1.1), vanishes to infinite
order in the coordinate $a_{m+1, \beta}$.
\medskip

(2)  If $f(0) \neq 0$, then there is a solution $u$ to (1.1) whose
Gevrey order is exactly $(m+1) \big/ (s+1),$ in the 
coordinate $a_{m+1, \beta}$ at the origin.
\endproclaim
\medskip

\remark{Remarks}  We now see clearly that the non-analyticity is an
interplay of the difference between the non-symplectic stratum and the center,
i.e. $s < m$.  On the Heisenberg group, $m= s = 1$, that is there are no
non-symplectic strata.  In
[BG], $m=1, s=0$ and so on.
\endremark 
\medskip

It is crucial to point out that we must complexify $\lambda_{s+1,
\alpha}$ or else we may not have any solutions to the Schr\"odinger
equation with the prescribed bound if we restrict to real values of $\lambda_{s+1, \alpha}$.
However, $\lambda_{m+1,\beta}$ is picked real.  Complex values of
$\lambda_{s+1, \alpha}$ off course do not contribute to the Plancherel
formula for $N$, and it is only the real values of $\lambda \in \frak
n^*$ that plays a role in the $C^\infty$ theory. The observation that
one needs to complexify  parameters of the representation to understand
analytic hypoellipticity was also made by B. Helffer[H].

\demo{Proof of Theorem 2.3}  Let $\rho > 0$.  We recall that there is
a natural dilation $\delta (\rho)$ on the graded
Lie Algebra $\frak n$, [F] and [RS].  Since $\frak n^*$ is dual to
$\frak n$, it inherits the same gradation as $\frak n$, and thus 
$$
\delta(\rho)\lambda = (0, \dots, 0, \rho^{s+1} \lambda_{s+1,\alpha } ,
0, \dots, 0, \rho^{m+1 } \lambda_{m+1, \beta }).\tag 2.21 
$$
Next, we have for $X_j \in \Cal B, a \in N, t\in \Bbb R^n,$
$$
X_j \pi_\lambda (a) f (t) = \pi_\lambda (a) d\pi_\lambda (X_j) f
(t). \tag 2.22
$$
The left side of (2.22) is 
$$
\align
 & \frac{\partial} {\partial \tau}\pi_\lambda (a e^{\tau X_j}) f
\big|_{\tau = 0}\\
 &=  \pi_\lambda (a) \frac{d}{d\tau} \pi_\lambda (e^{\tau
X_j})f\big|_{\tau = 0} = \pi_\lambda (a) d\pi_\lambda(X_j) f.
\endalign
$$
Thus (2.22) is established.  Define for $a\in N$, 
$$ u(a) = \int^\infty_0 \left[\pi_{\delta(\rho)\lambda} (a) f_\rho (t)
\right]\big|_{t=0}  e^{-M\rho^{s+1}} d\rho \tag 2.23
$$
where $f_\rho(t) = f(\rho t)$, and $M> 0$ is a suitable large number.  From
(2.15), the definition of $\delta(\rho)$ and the hypothesis (2.20), we
see that if $a\in U, U $ a suitably picked neighborhood of the origin
in $N$, and $M$ is picked large enough, the integral in (2.23) will
converge absolutely.  Now by (2.22),

$$
\sum^p_{j=1} X_j^2 u = 
\int^\infty_0 \left[\pi_{\delta(\rho)\lambda} (a) d\pi_{\delta(\rho)
\lambda } \left(\sum^p_{j=1} X_j^2\right) f_\rho(t)\right]\biggr|_{t=0} e^{-M\rho^{s+1}} d\rho.
$$
But, 
$$
d\pi_{\delta(\rho)\lambda} \left(\sum^p_{j=1} X_j^2\right) f_\rho
\equiv 0. \tag 2.24
$$
This is because,
  
$$
\align
d\pi_{\delta(\rho)\lambda}\left(\sum^p_{j=1} X_j^2\right) f_\rho =&
\sum^n_{j=1} \left(\frac{\partial}{\partial t_j} + i \left(
\rho^{s+1}\lambda_{s+1, \alpha } \tilde p_j (t) + \rho^{m+1}
\lambda_{m+1, \beta } \tilde q_j (t)\right)\right)^2f_\rho\\ 
 -& \sum^r_{j=1} \left(\rho^{s+1}\lambda_{s +1, \alpha } p_j (t) +
\rho^{m+1}\lambda_{m+1, \beta } q_j (t)\right)^2 f_\rho.
\endalign
$$
Since $\tilde p_j, p_j$ are homogeneous of degree $s$, and $\tilde
q_j$ and $q_j$ homogeneous of degree $m$, by changing variables $w =
\rho t$, we may write the expression above as, $\rho^2 A f(w)$.  By
hypothesis this vanishes.  Thus, $u(a)$ is a solution.  We now denote
the variable $a_{m+1, \beta}$ on $N$, by simply $a_{m+1}$.

We calculate,
$$
\frac{\partial^\sigma u }{\partial a^\sigma_{m+1}} ( 0, \dots, 0, a_{m+1})\big|_{a=0}.
$$
Let $\tilde a = ( 0, 0, \dots, 0, a_{m+1, \beta}).$
Then we claim,
$$
\pi_{\delta(\rho)\lambda} (\tilde a ) f_\rho (t) =
 e^{i\lambda_{m+1, \beta } \rho^{m+1}a_{m+1}} f_\rho (t).\tag 2.25
$$
Since $\tilde a\in V_{m+1}$, it commutes with every element in $N$.
Thus, if $t= (t_1, \dots, t_n),$
$$
\align
(t, 0) \bullet \tilde a =& (t, 0, \dots, 0, a_{m+1, \beta})\\
  =& (0, 0, \dots, 0, a_{m+1, \beta}) \bullet (t, 0).
\endalign
$$
Using (2.10), (2.25) follows.

Thus, combining (2.23) and (2.25),
$$
u (0, a_{m+1}) = \int^\infty_0  e^{i\lambda_{m+1, \beta} \rho^{m+1}
a_{m+1}} f(0) e^{-M\rho^{s+1}} d\rho.\tag 2.26
$$
Note now that it is crucial $\lambda_{m+1, \beta } \in \Bbb R$ for the
integral to converge.  If $f(0) = 0$, then (2.26) vanishes with all
derivatives.  If $f(0) \neq 0$, then 
$$
\frac{\partial^\sigma u}{\partial
a^\sigma_{m+1}}(0, a_{m+1})\big|_{a=0}= f(0) \lambda^\sigma_{m+1, \beta}
i^\sigma\int^\infty_0 \rho^{(m+1)\sigma} e^{-M\rho^{s+1}} d\rho.
$$
Changing variables we see easily that for $C>0$,
$$
\bigg| \frac{\partial^\sigma u }{\partial a^\sigma_{m+1}} 
(0, a_{m+1}) \bigg| \sim C^\sigma (\sigma!)^{(m+1)/(s +1)}.
$$
We also used the fact that $\lambda_{m+1,\beta}\neq 0$ in the last step.
\enddemo

$\hfill \blacksquare$

The next two propositions seek to establish a link between non-symplectic
strata in the characteristic set $\Sigma$ of $L$ and the appearance
of potential wells in the Schr\"odinger equations attached to the
derived representations $d\pi_\lambda(L)$. We have only been able to establish
this link when the characteristic set itself is non-symplectic as in
the situation of [BG] and that too in a stronger situation. The difficulty
in establishing this link when the non-symplectic stratum is immersed
as in [HH2] is in the combinatorial aspects of the Campbell-Hausdorff
formula.

\proclaim{Proposition 2.4} Assume that the vector space $V_1$ contains
an element from the center of $\frak n$. Then the characteristic set
$\Sigma$ of $L$ is non-symplectic. The converse statement does not hold.
\endproclaim

\demo{Proof} Since symplectic properties of $\Sigma$ remain invariant
under linear combinations, we may suppose that in our basis $\Cal B$
for $V_1$, the element $X_p$ is in the center of $\frak n$. Denote
the symbols of the vector fields $X_i$ in the basis $\Cal B$ by
$\sigma (X_i)$. By hypothesis $[X_i, X_p]=0$ for all $X_i\in \Cal B$.
Thus, $\{\sigma(X_i),\sigma(X_p)\}=0$, for every $X_i\in \Cal B$. Here
$\{,\}$ denotes the Poisson bracket.
Since $\Sigma$ is defined by the set of equations $\sigma(X_i)=0,\ X_i
\in \Cal B$, we conclude right away that the presence of the element
$X_p$ from the center of $\frak n$ forces $\Sigma$ to be non-symplectic.

Now on the other hand consider 
$$
X_1=\frac{\partial}{\partial x_1},\, X_2=\frac{\partial}{\partial x_2}+
x_1\frac{\partial}{\partial x_3},\, X_3=\frac{\partial}{\partial x_4}+x_1
\frac{\partial}{\partial x_5}
$$
The characteristic set of the operator $L=\sum^3_{i=1}X_i^2$ is 
non-symplectic. However in the Lie algebra generated
by the $X_i$, the vector space $V_1$ does not contain an element from the
center. Note by considering a solution of the form $u(x_1, x_2,x_5)$
the situation reduces to the example in [BG] and thus the operator
$L$ is not analytic hypoelliptic.

\enddemo
$\hfill\blacksquare$

\proclaim{Theorem 2.5} Assume that the basis $\Cal B$ contains an element
from the center of $\frak n$. Furthermore set $\lambda_{s+1,\alpha}=\mu$
and $\lambda_{m+1,\beta}=1$. Then we have

$$ 
d\pi_\lambda (L) = \sum^n_{j=1}\left(\frac{\partial}{\partial t_j} +
i 
\tilde q_j (t)\right)^2
 - \sum^{r-1}_{k=1} q^2_k(t) -\mu^2.
$$
\endproclaim

\demo{Proof} Let us assume that the basis element $X_p\in \Cal B$ is
an element of the center of $\frak n$. We select $\lambda_1$, such that,
$\lambda_1(X_p)=\lambda_{s+1,\alpha}=\lambda_{1,p}=\mu$ and $\lambda_1(X_i)
=0,\, i\neq p$ where $X_i\in \Cal B$. Further let 
$\lambda_1|_{V_i}=0, i\neq 1$. Thus we are in the 
situation $s=0$. Since $s=0$, by construction $\Cal S_1=\phi$. 
Now $X_p\notin \Cal S_2$,
because $\langle \lambda_2,[Z,X_p]\rangle=0$ for any $Z\in \frak n$.
Thus $X_p\in \Cal B\backslash \Cal S$ and thus $X_p\in W$. We
conclude that $\lambda_1|_W\neq 0$. This means that $\langle
\lambda_1,\frak h
\rangle \neq 0$. Thus $\langle \lambda_1,\frak h \rangle$ contributes
to the exponential in (2.10) and thus to (2.15) which is the 
induced representation $\pi_\lambda$. 
From (2.12)
we conclude that (2.13) takes the form
$$ \text{(component)}_{1,p}=a_p. \tag 2.27$$
Next we focus on (2.14). Since $X_p$ lies in the center once again
we see right away from the Campbell-Hausdorff formula that the second
sum on the right in (2.14) will contain no linear terms involving $a_p$. 
That is (2.14) in the current situation takes the form
$$
\align
\text{(component)}_{m+1,\beta} =& \sum^n_{j=1} \tilde q_j (t) a_j +
\sum^{r-1}_{j=1} q_j(t) a_{j+n}\\
 +& \sum_{j\geqslant p + 1} b_j' (t) a_j + \sum q_{j_1j_2\dots
j_k} (t) a_{j_1} a_{j_2} \dots a_{j_k}.\tag 2.28
\endalign
$$
Recall that $p=r+n$, and the second sum on the right above stops at 
$r-1$ because $X_p$ commutes with every element of $\frak n$.
Using now (2.27) and (2.28) instead of (2.17) and (2.18) and setting,
$\lambda_{m+1,\beta}=1,\, \lambda_{1,p}=\mu$ we easily get our conclusion.
\enddemo
$\hfill\blacksquare$

\remark{Remark 2.6}  Though the results of this section were derived
under the assumption
$$
\lambda_2 = (0, 0 \dots, 0, \lambda_{m+1, \beta })
$$
an examination of the proofs shows that one could take for $\lambda_2
\in V^*_{m+1},$
$$
\lambda_2 = (0,\dots, 0, \lambda_{m+1, \beta_1 }, \lambda_{m+1,
\beta_2 } , \dots, \lambda_{m+1, \beta_k }, 0 )
$$
provided, (2.4) holds.  By allowing this flexibility it is technically
easier to verify (2.20).  This remark will become clear when we
consider examples in Section 4.
\endremark

\subhead \S 3.  Bounds for the Schr\"odinger Equation \endsubhead
\medskip

We now try to prove the bound given in (2.20) for the spectral problem
in Theorem (2.3).  We begin with some remarks.
\medskip

When $n=1$, the ODE case, we observed by
Proposition (2.2), that the spectral problem is, for $a_j, b_j \in
\Bbb R$,
$$
\frac{d^2 f}{dt^2} - \sum^r_{j=1} \left(\lambda a_j t^s + \mu b_j
t^m\right)^2 f = 0.
$$
with $s< m$.  We pick $\mu =-1$, and we see easily that we get the
problem for $\lambda \in \Bbb C, a_j, b_j \in \Bbb R,$
$$
-\frac{d^2 f}{dt^2} + \sum^r_{j=1} \left(b_j t^m - \lambda a_j
t^s\right)^2 f = 0.
$$
This problem has been the principal problem treated in the literature [BG],
[C1, 2], [HH1, 2], [Hos], [LR], [H], [CC], [O].  Thus under (2.3) and
(2.4) and when $n=1$, the problem of analytic hypoellipticity on a
nilpotent group is essentially settled, though all authors treated
$r=1$.  The main point is when $n=1$, all potentials of the
Schr\"odinger equation are confining.  For $n\geqslant 2$, we get for the
first time non-confining potentials and the problem is far more
difficult even though we are asking for weak bounds (2.20).  We have
been only able to treat special cases of the Schr\"odinger equation. At
the end of this section we present a theorem of F. Treves that hints
that (2.20) is true in general.  Unfortunately we have been unable to
capitalize on Treves's theorem.
\medskip

We assume that (2.3), (2.4) holds and the vector fields in $\Cal S$
commute.  Then by Section 2, we arrive at the Schr\"odinger equation,
(after setting $\lambda_{m+1, \beta } = - 1, \lambda_{s+1, \alpha} =
\sigma)$ 
$$
-\Delta f + \sum^r_{j=1} \left(\sigma p_j (t) - q_j (t)\right)^2 f(t)
= 0.
$$
Let us assume that for all pairs $q_j (t), p_j (t)$ if $q_j(t) \neq 0$,
then $p_j (t) \equiv 0,$ and if $p_j(t) \neq 0$, then $q_j \equiv 0$.
Then the equation above becomes,
$$
-\Delta f + \left(\sum^{r'}_{j=1}q^2_j(t) + \sigma^2\sum^{r''}_{j=1}
p^2_j (t)\right) f=0.
$$
Notice for $\sigma \in \Bbb R$, we need not have a solution with the
prescribed growth (2.20) and
we are forced to pick $\sigma \in \Bbb C$.  In fact we take $\sigma =
i\sqrt{\lambda},  \lambda > 0$ and we get the problem, $t\in \Bbb R^n$
$$
-\Delta f + (q(t) - \lambda p(t))f = 0 , \lambda >0,
$$
$q\geqslant 0$, homogeneous polynomial of degree $2m, p(t) \geqslant
0$, homogeneous polynomial of degree $2s$, $m>s$, and both polynomials
having real coefficients.  Unlike the ODE case $n=1$, we cannot
anymore expect $f(t)$ to decay.
\medskip
\noindent{\bf Example (3.1)}  Let $h(t)$ be a homogeneous harmonic
polynomial in $\Bbb R^n$ of degree $k$.  Let, 
$$
f(t) = e^{-h^2/2}.
$$
then $f(t)$ is a solution to,
$$
-\Delta f + (q(t) - p (t)) f (t) =0,
$$
$q(t) = h^2 |\nabla h |^2, p(t) = |\nabla h|^2,$ thus $\deg p =
2(k-1), \deg q = 2 (2k-1)$.  In $\Bbb R^2$, if we select $h(t) = t^2_1
- t^2_2$, we notice that along the ray $t_1 = t_2, f (t) \equiv 1$.
Thus decay in the eigenfunction problem above is not possible and
(2.20) does not demand it either.
\medskip

We have,

\proclaim{Theorem 3.2}  Let us consider in all $\Bbb R^n$,
$$
-\Delta f + (q (t) - \lambda p (t)) f(t) = 0,
$$
$q(t) \geqslant 0$ is a homogeneous polynomial of degree $2m, p(t)
\geqslant 0$ is a non-trivial, homogeneous polynomial of degree $2s,
m>s$.  We have the following.  If 
\medskip

(a)  $q(t)$ does not vanish on the unit sphere $S^{n-1}$,

or 
\medskip

(b)  $q(t)$ vanishes on $S^{n-1}$ and the maximal order of vanishing
of $q(t)$ restricted to $S^{n-1}$ is $j$, and
$$
\frac{2s - 2m + j}{j} + s < 0. \tag 3.1
$$
Then there exists $\lambda > 0 $, and a solution $f(t)$ to the
Schr\"odinger equation above,
such that,
\medskip

(1)  $f(t) > 0 , f(0) =1,$
\medskip

(2)  $0< f(t) < C ( 1+ |t|)^{(m+1)(n-1)}, t\in \Bbb R^n.$
\endproclaim
\medskip

The case $j =0$ is included in (3.1), and under assumption (a) we can
even conclude that $f(t)$ decays.
\medskip

The proof relies very strongly on Harnack's inequality.  In fact we
need a scale invariant version.

\proclaim{Lemma 3.3}  Let $g\geqslant 0.$  Let $g$ satisfy in $\Bbb
R^n$,
$$
-\Delta g + V (t) g =0.
$$
Let $B_R(t_0)$ be a ball centered at $t_0$ of radius $R$.

Assume, 
$$
R \sup_{B_{2R}(t_0)} | V(t)|^{1/2} \leq C_1. \tag 3.2
$$
Then,
$$
\sup_{B_R(t_0)} \leq C(C_1, n) \inf_{B_R(t_0)} g.
$$
\endproclaim
\demo{Proof}  This is a special case of Theorem 8.20 in
Gilbarg-Trudinger [GT].  Set $w(t) = g (Rt + t_0)$.

Then, $$
-\Delta  w + R^2 V(Rt + t_0) w(t) = 0 
$$
$t\in B_2 (0)$, where $B_2(0)$ is the ball of radius $2$ around the
origin.  Using (3.2) and the result in [GT], our conclusion follows.
\enddemo

$\hfill \blacksquare$

We will display the proof of theorem (3.2) for $n=2$, and indicate the
changes necessary for the higher dimensional version.  We do so to
keep the ideas transparent, the idea being the same in all dimensions.
Perhaps this will help to eventually get rid of assumption (3.1) which
is our goal.  The proof in higher dimension involves
notational complications and use of Whitney's lemma.  
\medskip

We begin with an approximate problem.  Let $\phi_N$ be a smooth,
radial cut-off function for $N$ large, given by,
$$\phi_N(t) =  \cases 1, & |t| < N\\ 0, & |t| > 2N.\endcases.
$$
Fix, $\lambda \in \Bbb R$.  We will also need polar coordinates $(r, \theta)$.  Then $q(t) = r^{2m} q_1 (\theta),
p(t) = r^{2s} p_1 (\theta)$.

Define
$$
T= \{ t : q(t) - \lambda p (t) < 1, |t|\geqslant 100\}.\tag 3.3
$$
The set $T = \overset{\ell}\to{\underset{\sigma = 1}\to \cup}
T_\sigma$, where each $T_\sigma$ is a tapering tube, that is, there is
a central axis, and transverse to this axis the diameter of $T_\sigma$
falls off.  Assume that we have rotated axes, so that $q_1(0) = 0$ and
let this correspond to $T_1$.  Let $k$ be the vanishing order of $q_1$ at
$\theta = 0$.  Then the set $T_1$ is bounded by the curves for large
$r$, given by,
$$
r^k |\theta|^k = c(\lambda r^{2s-2m+k} + r^{k-2m}).
$$
Taking, $k-th$ roots, we have the bounding curves can be taken to be,
for large $r$,
$$
t_2 = \pm c(1+|\lambda|)^{1/k} r^{(2s-2m + k)/k}, t=(t_1, t_2). \tag 3.4
$$
By (3.1), $(2s - 2m + k)/k < 0$.
\medskip

We now consider an approximate problem,
$$
-\Delta f_N + (q(t) -\lambda_N p(t) \phi_N (t)) f_N (t) = 0.\tag 3.5
$$
The new potential $q(t) - \lambda_N p (t) \phi_N(t)$ is bounded below
so we hope to find a non-negative solution to (3.5) by a variational
approach.
\proclaim{Lemma 3.4}  Let the weaker condition $2s-2m + j< 0$ hold.
Then there exists $\lambda_N > 0, f_N > 0$, with 
$$
\int_{\Bbb R^n} |\nabla f_N |^2 + \int_{\Bbb R^n} |f_N|^2 < + \infty
\tag 3.6
$$
and such that $f_N$, satisfies (3.5).
\endproclaim

\demo{Proof}  We will see that $f_N$ satisfies a stronger condition
than (3.6).  Fix a value $\lambda$, and consider the sets
$$
\align
S_1 & = \{ t: q(t) - \lambda p (t) \phi_N(t) \geqslant 1\}\\
S_2 & = \{ t: q(t) - \lambda p (t) \phi_N (t) < 1\}.
\endalign
$$

We shall drop the subscript $N$ in $f_N$ from now on.
\medskip

Now,
$$
\align
-& \Delta f + (q - \lambda p \phi_N) f= - \Delta f + (q-\lambda
p \phi_N)(\chi_{S_1} + \chi_{S_2}) f\\
=& -\Delta f + [(q-\lambda p \phi_N) \chi_{S_1} + \chi_{S_2} ] f + [(q
- \lambda p \phi_N)\chi_{S_2} - \chi_{S_2} ] f\\
=& - \Delta f + V_1 (t) f - V_2 (t) f,\tag 3.7
\endalign
$$

where,

$$
\align
 V_1 (t) & = (q - \lambda p \phi_N) \chi_{S_1} + \chi_{S_2} \geqslant
1,\\
V_2 (t) & = \chi_{S_2} - (q-\lambda p \phi_N) \chi_{S_2} \geqslant
0.\tag 3.8 
\endalign
$$

Moreover,
$$
\| V_2\|_{L^\infty} \leqslant c (\lambda, N)\tag 3.9
$$

Define the Hilbert space $\Cal H$, to be the completion of the vector
space $W$,

$$
W= \{\phi \in C^\infty_0 (\Bbb R^n) : \int_{\Bbb R^n} |\nabla\phi |^2
+ \int_{\Bbb R^n} |\phi|^2 V_1 < + \infty\},
$$
with norm
$$
\| \phi\|^2_{\Cal H} = \int_{\Bbb R^n} |\nabla \phi |^2 + |\phi|^2 V_1.
$$
Let $\Cal H^*$ be the dual to $\Cal H$.  Then, $A= (-\Delta +
V_1)^{-1}$ is a bounded operator from $\Cal H^* \to \Cal H$, and in
fact $\| A \|\leqslant 1$.  We show that $A(V_2 f) = Kf$ is a compact
operator on $\Cal H$.  To demonstrate this, note by (3.9), and since
$V_1\equiv 1$ on the support of $V_2$,
$$
K_1 : \Cal H \to \Cal H^* , K_1 f = V_2 f,
$$
is a bounded operator.  Now $S_2\subseteq T,$ for $r\geqslant 100$.
Thus $S_2$ is bounded by the curves (3.4) which are tapering in the
$t_2$ direction.  This tapering will produce the desired
compactness.  In the general $\Bbb R^n$ case it is seen that the
components of the set $T$, tapers at least in one distinguished
direction and that is all that is required in the ensuing proof.  We
now show $K_1$ is a compact operator.  For $\varepsilon > 0$, we
prove, there exists $R_\varepsilon$, such that for $C>0, C = C
(\lambda, N)$,
$$
\int_{|t| > R_\varepsilon} |f|^2 V_2  < C \varepsilon \| f \|^2_{\Cal
H}\tag 3.10
$$
Once (3.10) is established, then an application of Rellich's lemma on
the compact part $\{|t| < R_\varepsilon\}$ and a diagonal process
establishes the compactness.  Since the support of $V_2$ is contained
in $T$, we have from (3.9),
$$
\int_{|t|> R_\varepsilon} |f|^2 V_2 \leqslant \sum^\ell_{\sigma = 1}
\int_{\{|t| > R_\varepsilon\} \cap T_\sigma} |f|^2.\tag 3.11  
 $$
So its enough to establish (3.10) for $\sigma = 1$.  Let $\beta = (2m
-2 s -k)/k > 0$.  Pick $\alpha > 0, \alpha < 1/2$ and $2\alpha <
\beta$.  Let $\psi \in C^\infty_0 (\Bbb R)$,
$$
\psi (\tau) = \cases  1, & |\tau | \leqslant 2\\
 0, & |\tau|\geqslant 4.\endcases
$$
Note first that, because $\alpha < \beta$,

$$
B_1 \cap \{ |t| > C_1(\lambda)\} \subseteq \{(t_1, t_2) : |t_2|
|t_1|^\alpha \leq 1\}. \tag 3.12
$$
This follows because the set on the left in (3.12) by virtue of (3.4)
is bounded by the curves, $t_2 = \pm c ( 1+ |\lambda |)^{1/k} | t_1
|^{-\beta}$.  For $|t|\geqslant c_1(\lambda)$, we have for $\alpha<
\beta, c ( 1 + |\lambda |)^{1/k} |t_1|^{-\beta} \leqslant
|t_1|^{-\alpha}$.  (3.12) follows.  Going back to the right side of
(3.11), we have, by (3.12),
$$
\int_{\{|t| > R_\varepsilon\} \cap T_1} |f|^2 \leqslant \int_{\{|t| >
R_\varepsilon\} \cap T_1} |f\psi(|t_2| |t_1|^\alpha) |^2 dt_1 dt_2.
$$
For fixed $t_1$, we apply Poincare's inequality to the $t_2$ variable
to the integral on the right above.  We see that,
$$
\align
\int_{\{ |t|> R_\varepsilon\} \cap T_1} |f\psi|^2 dt_2 &\leq c
R_{\varepsilon}^{-2\alpha} \int_{\{ |t|> R_\varepsilon\}} |\nabla_{t_2}
f|^2 \psi^2 dt_2 \\
&+ c\int_{\{|t| >R_\varepsilon\} \cap \{ 2< |t_2|
|t_1|^\alpha < 4\}} |\psi'|^2 |f|^2 |t_1|^{2\alpha} dt_2.
\tag 3.13\endalign
$$
Now, for $\alpha'$, such that $\alpha < \alpha'< 1/2$ we have for
$R_\varepsilon$ large,
$$
|t_1|^{2\alpha} < \varepsilon |t_1|^{2\alpha'}.\tag 3.14
$$
We now show that if $2< |t_2| |t_1|^\alpha < 4,$
$$
|t_1|^{2\alpha'} \leq q(t) -\lambda p(t). \tag 3.15
$$
\medskip

{\bf Case 1: }  If $s=0, p(t) = c$.  Now when $2< |t_2| |t_1|^\alpha <
4, q (t) \sim |t_2|^k |t_1|^{2m-k} \sim |t_1|^{2m-(1+\alpha) k}$.
Thus we need to have, $2m -(1 +\alpha) k > 2 \alpha'$.  Since $k\neq
2m,$ by choosing $\alpha, \alpha'$ sufficiently small and positive we
easily verify (3.15).
\medskip

{\bf Case 2: }  If $s\neq 0$, then we have for $\alpha'$ small,
$2\alpha' < 2s$, and thus on $2< |t_2| |t_1|^\alpha < 4, |t| \geq
R_0(\lambda),$
$$
|t_1|^{2\alpha'}+ \lambda p (t) \leq c |t_1|^{2s} < |t_1|^{2m-k}
|t_2|^k \sim q(t).
$$
The last inequality holds whenever $|t_2| |t_1|^\beta \geqslant c$.
Since $|t_1|\geqslant 10$, and $|t_2| |t_1|^\alpha \geqslant c, |t_2| |
t_1|^\beta \geqslant c.$

Thus (3.15) holds in all cases.  Combining (3.14) and (3.15), for $|t|
\geqslant R_\varepsilon$,
$$
|t_1|^{2\alpha} \leq \varepsilon (q(t) - \lambda p (t) \phi_N (t)). \tag 3.16
$$
Thus for large $R_\varepsilon$, in (3.13),
$$
\int_{\{|t| > R_\varepsilon\} \cap T_1} |f\psi |^2 dt_2 \leqslant c
\varepsilon \left(\int_{\{|t_1|^\alpha |t_2| \leq 4\}}(|\nabla f|^2 +
|f|^2 V_1)dt_2\right).
$$
Integrating both sides above in the $t_1$ variable, we easily get
(3.10).  Thus the self-adjoint, compact operator, $(-\Delta +
V_1)^{-1} (V_2f) = Kf$, has a largest positive eigenvalue
$\mu_1^{-1}$, and eigenfunction $f\in \Cal H$, such that
$$
Kf = \mu^{-1}_1 f.
$$
We re-write this equation as, 
$$
-\Delta + V_1 f - \mu_1 V_2 f =0.\tag 3.17
$$
We now choose $\lambda$, so that $\mu_1 = 1$, and then with $\mu_1 =
1$, from (3.7) we see that $f$ will satisfy our Schr\"odinger equation.
The value of $\lambda$, which gives $\mu_1 = 1$, is the desired
$\lambda_N$.  Since $f\in \Cal H$, and $V_1 \geqslant 1$, we also see,
$$
\int_{\Bbb R^2} (|\nabla f|^2 + |f|^2) \leq \| f \|^2_{\Cal H} <
+\infty.
$$
To see we can select $\lambda_N$, so that $\mu_1 = 1$, we write
(3.17) in a variational framework.
$$
\mu_1 = \inf_{\phi\in C^\infty_0(\Bbb R^2)} \frac{\int_{\Bbb R^2}
|\nabla \phi|^2 + V_1 \phi^2}{\int_{\Bbb R^2}\phi^2 V_2}. \tag 3.18
$$
We claim:

(a)  If $\lambda \to +\infty, \mu_1 \to 0$.

(b)  If $\lambda \to - \infty, \mu_1 \to + \infty$.

Thus from (a), (b), there exists a $\lambda_N$, for which $\mu_1 = 1$.
This will be necessarily positive.  This is because, for the
eigenfunction $f$ and the $\lambda_N$, we must have,
$$
-\Delta f + q(t) f - \lambda_N p (t) \phi_N f = 0.
$$
Thus,
$$
\lambda_N = \frac{\int_{\Bbb R^2}|\nabla f|^2 + q f^2}{\int_{\Bbb R^2}
p (t) \phi_N f^2}.
$$
Thus $\lambda_N > 0$, since $q(t), p(t), \phi_N(t) \geqslant 0$, with
$p(t)\not\equiv 0$.  It also follows from (3.18) that $f>0$.  We now
check our claim.

{\bf Case (a): }  Since,
$$
V_2 (t) = \chi_{S_2} - (q - \lambda p \phi_N) \chi_{S_2},
$$
as $\lambda \to + \infty$, the support of $V_2$ will contain a fixed ball
$B$ of radius $r_0$, such that $V_2\geqslant M$ on this ball.  The
ball $B$ can be fixed since, for $\lambda_1 <
\lambda_2$,
$$
\{t : q(t) - \lambda_1 p(t) \phi_N (t) < 1\} \subseteq \{ t:
q(t)-\lambda_2 p(t) \phi_N (t) < 1\}.
$$
On this ball $V_2\to + \infty$ as $\lambda \to \infty$.  Pick $\phi
\in C^\infty_0(B)$.  Then substituting this test function into (3.18),
we notice on $B, V_1\equiv 1,$ for large $\lambda$, from the
definition of $V_1$, while the denominator goes to positive infinity,
thus $\mu_1 \to 0$.
\medskip

{\bf Case (b): }  Let $\lambda \to - \infty$.  Now notice for $\lambda
\leqslant 0, q - \lambda p = q + |\lambda | p \geqslant 0$, and thus
$\| V_2\|_\infty \leq 2$.  Secondly by examining the proof that led to
(3.15), for $\lambda \leqslant 0$, we see easily we have for
$\varepsilon > 0$, there exist $R_0$, independent of $\lambda$, such
that, for $\lambda \leqslant 0$ and $2 < |t_2 | |t_1|^\alpha < 4$ and
$|t|> R_0$,
$$
|t_1|^{2\alpha} \leqslant \varepsilon ( q (t) + |\lambda | p (t)
\phi_N (t)) = \varepsilon V_1.
$$
Thus for any $\phi \in C^\infty_0$, we see if $\lambda \leqslant 0$,
$$
\int_{\{ |t| > R_0\}} |\phi |^2 V_2 \leq C\varepsilon \| \phi\|^2_{\Cal H} \tag 3.19
$$
with $C$ independent of $\lambda$ and even $N$.  Now for $|t| < R_0$,
as $\lambda \to -\infty$,
$$
|S_2 \cap \{ |t| < R_0\} |\to 0.\tag 3.20
$$  
Thus,
$$
\align
\int_{\Bbb R^2} |\phi|^2 V_2 =& \int_{T\cap\{|t|> R_0\}}+ \int_{T\cap \{
|t|\leq R_0\}} | \phi|^2
V_2 \\
    \leq&  c\varepsilon \|\phi\|^2_{\Cal H} + \left(\int_{T\cap\{ |t|
\leq R_0\}} V_2^2\right)^{1/2}\left(\int_{\{|t|\leqslant
R_0\}}|\phi|^4\right)^{1/2}.
\endalign
$$
From (3.20), 
$$
\align
\leq &c\varepsilon \|\phi\|^2_{\Cal H} + c \varepsilon
\left[\int_{\{|t|\leq R_0\}}(|\nabla \phi|^2 + \phi^2)\right]\\
  \leq &c\varepsilon \|\phi\|^2_{\Cal H}.
\endalign
$$
Thus, for any $\phi \in C^\infty_0 (\Bbb R^2)$, as $\lambda \to -
\infty$, with a uniform $c$ we have,
$$
\frac{c}{\varepsilon} \leqslant \frac{\|\phi\|^2_{\Cal H}}{\int_{\Bbb
R^2}|\phi|^2 V_2}.
$$
Claim (b) follows and our lemma is proved.
\enddemo

$\hfill \blacksquare$

\proclaim{Lemma 3.5}  Under (3.1), there exists a finite, positive
constant $C>0$, independent of $N$, such that the eigenvalues
$\lambda_N$ of Lemma (3.4) satisfy,
$$
C^{-1} < \lambda_N < C\tag 3.21
$$
\endproclaim
\demo{Proof}  The upper bound in (3.21) is easy to come by.
Re-writing (3.18) with $\mu_1=1$, as we have done earlier, we see,
$$
\lambda_N = \inf_{\phi \in C^\infty_0(\Bbb R^2)}\frac{\int_{\Bbb R^2}
(|\nabla \phi |^2 + q\phi^2)}{\int_{\Bbb R^2} p \phi_N \phi^2}.\tag 3.22
$$
Now choose,
$$
\phi (t) = \cases 1 & |t|\leqslant 1\\0, & |t|\geq 2\endcases.
$$
and insert in (3.22).  One sees that $\phi_N \equiv 1$, on the support of
$\phi$, and hence we get easily $\lambda_N \leq C$, with $C$
independent of $N$.
\medskip

We now establish the lower bound for $\lambda_N$.  We know from Lemma
(3.4), that for the eigenfunction $f_N\equiv f$, we have,
$$
\lambda_N = \frac{\int_{\Bbb R^2} |\nabla f |^2 + q f^2}{\int_{\Bbb
R^2} p(t) \phi_N f^2}\tag 3.23
$$
We will show, for $C>0$, independent of $N, f$
$$
\int_{\Bbb R^2} p (t) \phi_N f^2 \leq c\int_{\Bbb R^2} q f^2\tag 3.24
$$
Combining (3.23) and (3.24) the lower bound for $\lambda_N$ follows.
First we note, by Lemma (3.3) and by the fact just established that
$\lambda_N \leqslant C$,
$$
\int_{|t|\leq 100} p\phi_N f^2\leq c\sup_{|t|\leq 100} f^2\leq
c\inf_{|t|\leq 100} f^2 \leq c \int_{|t|\leq 100} f^2 q.\tag 3.25
$$
Thus we show (3.24) for $|t|\geqslant 100$.  We revert to polar
coordinates again.  Fix $r$, and a circle $|t|=r$.  We show for
$r\geq 100$.

$$\int^{2\pi}_0 p \phi_N f^2 d\theta \leq c\int^{2\pi}_0 q f^2 d
\theta,\tag 3.26
$$
with $c$ independent of $r, N, f$.  Multiplying (3.26) by $r$, and
integrating over $r\geq 100$, we get,
$$
\int_{|t|\geq 100} p\phi_N f^2 dt \leqslant c \int_{|t|\geq 100} qf^2
dt.\tag 3.27
$$
Combining (3.27) and (3.25), (3.24) follows.  We now use the notation
$r\theta = (r \cos \theta, r \sin \theta)$.  Using the homogeneity of
$p, q$, we write (3.26) as, 
$$
\phi_N(r) r^{2s}\int^{2\pi}_0 f^2(r\theta) p_1(\theta) d\theta \leq c
r^{2m} \int^{2\pi}_0 f^2(r\theta) q_1(\theta) d \theta.\tag 3.28
$$
For a fixed $r\geq 100$,
$$
\Gamma_r = \{\theta : a r^{2s} \geq r^{2m} q_1(\theta)\},\tag 3.29
$$
where $a=\|p_1\|_{L^\infty(S^1)}$.  Now $r\Gamma_r$ (the dilate of
$\Gamma_r$ by a factor of $r$) consists of a union of $\ell$ arcs,
each arc transversal to one of the tapering tubes $\tau_\sigma$, defined
earlier of which there are $\ell$ in number.  Let us set $r\Gamma_r
\cap T_1 = G_{r,1}$.  Now note on the complement of $\Gamma_r$ in
$S^1, \Gamma^c_r$ we have,   
$$
\phi_N(r) r^{2s}\int_{\Gamma^c_r} f^2 (r\theta) p_1(\theta) d\theta
\leq r^{2m}\int_{\Gamma^c_r} f^2(r\theta) q_1(\theta) d\theta.
$$
So to prove (3.28), it is enough to show,
$$
\phi_N(r) r^{2s}\int_{\Gamma_{r,1}} f^2(r\theta) p_1(\theta) d\theta
\leq c r^{2m}\int^{2\pi}_0 f^2(r\theta) q_1(\theta) d\theta,
\tag 3.30
$$
where
$$
\Gamma_{r,1} = \{\theta : r\theta \in G_{r,1}\}.
$$
Now let $\widetilde \Gamma_{r,1}$, be an arc, such that
$|\widetilde \Gamma_{r,1}|=5|\Gamma_{r,1} |,$ and $\widetilde\Gamma_{r,1}$
disjoint from $\Gamma_{r,1}$ but having common end-point with
$\Gamma_{r,1}$.  Thus $\widetilde \Gamma_{r,1}$ is adjacent to
$\Gamma_{r,1}$ with five times the length of $\Gamma_{r,1}$.  Now the
arc $r\Gamma_{r,1}$ lies between the curves, $|t_2| = c|t_1|^{(2s-2m+k)/k},
c=c(a)$ (see also (3.4)) and follows from the definition (3.29).
Thus, ($|\bullet|$ denotes arc-length) 
$$
r|\Gamma_{r,1} | \approx r^{(2s-2m+k)/k} \approx r
|\widetilde\Gamma_{r,1}|.\tag 3.31
$$
Moreover using (3.29) again, we see on $\widetilde\Gamma_{r,1}$,
$$
ar^{2s}\leqslant r^{2m} q_1(\theta) \leqslant c a r^{2s}, \tag 3.32
$$
$c$ independent of $r, \theta$.  In fact for $(t_1, t_2)\in
\widetilde\Gamma_{r,1},$
$$ |t_1|^{(2s-2m+k)/k}\leqslant |t_2| \leqslant 5
|t_1|^{(2s-2m +k)/k}.
$$
Now fix a ball $B, B \supseteq (r\Gamma_{r,1})\cup (r\widetilde
\Gamma_{r,1})$ and diam $B = |r\Gamma_{r,1} \cup
(r\widetilde\Gamma_{r,1})|$.  Since we have (3.32), on $B,
|V|=|(q-\lambda_N p)|\leqslant c r^{2s}$.  By (3.31), it follows that
if $d=$ diameter of $B$, 
$$
d|V|^{1/2} \leqslant c r^{((2s-2m +k)/k+s)}\leq c,
$$
since (3.1) holds and $r\geqslant 100$.  Thus Lemma (3.3) applies.  We
have, for $r\geqslant 100$,
$$
\sup_{\Gamma_{r,1}} f^2(r \theta) \leqslant c \inf_{\widetilde\Gamma_{r,1}}
f^2 (r\theta). \tag 3.33.
$$
Going back to (3.30), using (3.33)
$$
\align
\phi_N (r) r^{2s}\int_{\Gamma_{r,1}} f^2(r\theta) p_1(\theta) d\theta
&\leqslant r^{2s} \inf_{\widetilde\Gamma_{r,1}} f^2 (r\theta)
\int_{\Gamma_{r,1}}p_1 (\theta) d\theta\\
 &\leqslant c r^{2s} |\Gamma_{r,1}|\inf_{\widetilde\Gamma_{r,1}}
f^2(r\theta).\tag 3.34 \endalign
$$
Thus the right side of (3.34) is bounded by
$$
\leqslant c r^{2s} \inf_{\widetilde\Gamma_{r,1}} f^2 (r\theta)
|\widetilde\Gamma_{r,1} |. \tag 3.35
$$
But on $\tilde \Gamma_{r,1}$ from (3.32) $r^{2s} \leqslant r^{2m}
q_1(\theta)$.  Hence (3.35) is bounded by,
$$
c\int_{\widetilde\Gamma_{r,1}} f^2 (r\theta) q_1 (\theta) r^{2m}
d\theta.
$$
This finishes the proof of (3.30).
\enddemo

$\hfill \blacksquare$

\proclaim{Lemma 3.6}  The following form of inequality (3.26) also
holds for $r\geqslant r_0, r_0$ independent of $N$, and for some
$\varepsilon  > 0$,
$$
r^{2s+\varepsilon} \int^{2\pi}_0 f(r\theta) d\theta \leqslant c \int^{
2\pi}_0 q_1(\theta) r^{2m} f(r\theta) d\theta.
$$
\endproclaim
\demo{Proof}  The proof is simply re-writing the proof of the previous
lemma.  In fact in the proof of the previous lemma we used a weaker
form of (3.1),
$$
\frac{2s-2m+j}{j} + s\leqslant 0.
$$
When $j=0$, since $q_1(\theta) > c > 0$, the corollary is obvious.
Thus assume $j>0$.  Define, for $\varepsilon >0$ (that will be picked
later),
$$
\Gamma_r=\{\theta : ar^{2s+\varepsilon} \geqslant r^{2m} q_1(\theta)
\}.
$$
The tubes, for example the tube $T_1$ is bounded by the curves $|t_2| =
|t_1|^{(2s+\varepsilon - 2m + k)/k}$.

Thus as before,
$$
r|\Gamma_{r,1} |\approx r^{(2s+\varepsilon - 2m + k)/k} \approx r
|\widetilde\Gamma_{r,1} |.
$$
Again on $\widetilde\Gamma_{r,1}$,
$$
a r^{2s+\varepsilon} \leqslant r^{2m} q_1(\theta) \leq c a
r^{2s+\varepsilon},
$$
$c$ independent of $\theta, r, N$.

Now fix a ball $B$ with diameter $|r\Gamma_{r,1} \cup
(r\widetilde\Gamma_{r,1})|$, and $B\supseteq r\Gamma_{r,1} \cup (r
\widetilde\Gamma_{r,1})$.  On $B$, the potential $V(t) = q(t) -
\lambda \phi_N (t) p (t)$, satisfies, $|V|^{1/2}\leq c
r^{s+\varepsilon/2}$.  If $d=$ diam $B$, then on $B, d|V|^{1/2}
\leqslant c r^\eta, \eta = ((2s - 2m +k)/k) +s + \varepsilon/2 +
\varepsilon/k.$  Choose $\varepsilon > 0$ so that $\eta \leqslant 0$,
which we can because of (3.1).  Thus we may apply Harnack's
inequality, Lemma (3.3) again to conclude for $f$, an identical
inequality as in (3.33).  Thus, 
$$
\align
r^{2s + \varepsilon} \int_{\Gamma_{r,1}} f(r\theta) d \theta &\leqslant
r^{2s+\varepsilon} \inf_{\widetilde\Gamma_{r,1}} f (r\theta)
|\widetilde\Gamma_{r,1} | \\
 &\leqslant \int_{\widetilde\Gamma_{r,1}} f(r\theta) q_1 (\theta) r^{2m}
d\theta,
\endalign
$$
because on $\widetilde \Gamma_{r,1}, r^{2m} q_1(\theta) \sim
r^{2s+\varepsilon}$.  
\enddemo

$\hfill \blacksquare$
\medskip

In $\Bbb R^n$, one proves this lemma by constructing the arcs
$\Gamma_{r,1}$ in the direction transversal to the direction of
tapering.  One direction always exists in which the diameter tapers
and we can apply the argument in that direction.

\medskip

Our $f_N \equiv f,$ satisfies $f_N > 0$, and solves, (3.5)
$$
-\Delta f_N + (q_1 (\theta) r^{2m} - \lambda_N p_1(\theta) r^{2s}
\phi_N (r)) f_N = 0.
$$
We may normalize $f_N$, such that $f_N(0) = 1$ and we assume this
from now on.  We integrate both sides of the equation above in
$\theta$.  We set, 
$$  
F(r) = \frac{1}{2\pi} \int^{2\pi}_0 f_N (r\theta) d\theta\tag 3.36
$$
By our normalization 
$$
F(0) = 1,\tag 3.37
$$
and we also see that,
$$
F'(0) = 0.\tag 3.38
$$
Further $F$ satisfies the ODE,

$$
\align
-(F'' + \frac{1}{r} F') 
&+ \left[\frac{r^{2m}\int^{2\pi}_0q_1(\theta)
f_N (r\theta) d\theta}
{\int^{2\pi}_0 f_N (r\theta) d\theta}\right.\\
      &-\left. \lambda_N \phi_N (r)
\frac{r^{2s}\int^{2\pi}_0p_1(\theta) f_N(r\theta) d\theta}
{\int^{2\pi}_0 f_N(r\theta) d\theta}\right] F(r) =0.
\endalign
$$
We set,
$$
\align
a(r) =& \frac{\int^{2\pi}_0 q_1(\theta)
f_N(r\theta)d\theta}{\int^{2\pi}_0 f_N(r\theta) d\theta}\\
b(r) =& \frac{\int^{2\pi}_0p_1(\theta) f_N(r\theta)
d\theta}{\int^{2\pi}_0 f_N (r\theta) d\theta}.\tag 3.39
\endalign
$$
Our ODE becomes
$$
-(F''+\frac{1}{r} F') + ( a(r) r^{2m} - \lambda_N\phi_N(r) \, b(r)
r^{2s}) F=0.\tag 3.40
$$
Now, by the Schwartz inequality and (3.36)
$$
\align
|F(r)|^2 &\leqslant c \int^{2\pi}_0 |f_N(r\theta) |^2 d\theta\\
|F'(r) |^2 &\leqslant c\int^{2\pi}_0 \big|\frac{\partial f_N}{\partial
r}\big|^2 d\theta.
\endalign
$$
Thus from (3.6) we have,
$$
\int^\infty_0 (|F'(r)|^2 + |F|^2) r dr < +\infty.\tag 3.41
$$
\proclaim{Lemma 3.7}  Let $F(r)$ on $[0, \infty)$ be defined by
(3.36).  Then $F(r)$ satisfies (3.37), (3.38), (3.39), (3.40) and
(3.41).  Moreover, for some $C > 0$,
$$
0 < F(r) \leq C,
$$
$C$ independent of $N$.
\endproclaim

\demo{Proof} Since $f_N> 0$, from (3.36) it follows that $F(r) > 0$.
That $F(r)$ satisfies (3.36), - (3.41) has been proved already, all
that remains is the upper bound for $F(r)$.  First for $c>0$,
independent of $N$ we have,
$$
\align
&0\leq b(r) \leq c , 0\leqslant a (r) \leqslant c,\\
&c_1 r^{2s+\varepsilon - 2m} \leq a (r) , r\geqslant R_0>0 \tag 3.42
\endalign$$
where $c_1>0, R_0$ are independent of $N$.  The first two inequalities
for $a(r),b(r)$
follows from the definitions (3.39), and the last inequality in (3.42) follows
from Lemma (3.6), once we re-write $a(r)$ as,
$$
a (r) = r^{2s-2m + \varepsilon} \frac{\int^{2\pi}_0 r^{2m} q_1 (\theta)
f(r\theta) d\theta}
{\int_0^{2\pi} r^{2s+\varepsilon} f(r\theta) d\theta}.
$$
Thus from (3.42), for $r\geqslant R_0$,
$$
r^{2m} a (r) \geq c_1 r^{2s+\varepsilon}. \tag 3.43
$$
Now using (3.42), (3.43), our potential satisfies for $r\geqslant
r_1$, ($r_1$ independent of $N$),
$$
c_1 r^{2s+\varepsilon} \leqslant a (r) r^{2m} - \lambda_N \phi_N (r)
b(r) r^{2s} \leqslant c_2 r^{2m}.
\tag 3.44
$$

We now show our solution $F(r)$ cannot have a local maximum for
$r>c_1$, with $c_1$ independent of $N$. Assume that a local maximum
is attained at $r=r_0>0$. Then,  $F'(r_0) =0,
F^{''}(r_0) \leqslant 0$.  Thus, from the ODE for $F$, since $F>0$,
$$
a(r) r^{2m} - \lambda_N \phi_N (r)  b(r) r^{2s} \leq 0,
$$
and using (3.42) and (3.43), we have,
$$
c_1r^{2s+\varepsilon} \leq \lambda_N r^{2s} \leq C r^{2s}.
$$

Thus $r_0\leqslant c_0$, $c_0$ independent of $N$. 
Fix $\epsilon<1/2$. We claim that $F(r)\leqslant\epsilon$ for
$r>r_1>c_1$. If there exists $r_2$ such that $F(r)>\epsilon$
for $r>r_2$ we contradict (3.41). Thus there exists $r_2$ such
that $F(r)\leqslant \epsilon$ for $r>r_2$, or there exists intervals
$[r_k,r_{k+1}]$, $r_k>c_1$ such that $F(r_k)=F(r_{k+1})=\epsilon$ and $F$ has
a local maximum in the interior of $[r_k,r_{k+1}]$ by Rolle's theorem.
If the second case occurs, we contradict the fact that $F(r)$ 
has no local maxima for $r>c_1$.
Our claim is proved. Since $F(0)=1$, and $F(r)\leqslant\epsilon<1/2$ for
$r>r_2$, either $F(r)\leqslant 1$ for all $r$ or $F(r)$ achieves its 
maximum at $r_0>0$, which is also a local maximum. But then this
$r_0\leqslant c_0$, $c_0$ independent of $N$. This proves our lemma.
Notice we have also proved that $F(r)\to 0$ as $r\to\infty$. This fact
can also be deduced by refined estimates in Davies[D].
 
\enddemo

$\hfill\blacksquare$

\proclaim{Lemma 3.8}  There exists a constant $C>0$, independent of
$N$, such that, for $t\in \Bbb R^2$,
$$
|f_N(t)| \leqslant C(1+|t|)^{(m+1)}.
$$
\endproclaim
\demo{Proof}  We use polar coordinates once again.  Define, for fixed
$r, r\geqslant 100$, 
$$
E_\theta = \{ \theta : f_N (r\theta) >  (1 + r)^{(m+1)}\}.
$$
The open set $E_\theta$ can be written as a disjoint union of
open arcs on $S^1$.  Thus, $E_\theta = \cup
I_\sigma$, where $I_\sigma$ are disjoint open arcs on $S^1$.  Now, for
$r\geqslant 100$, using Lemma (3.7),
$$
\align
|E_\theta| r^{m+1} &\leqslant \int_{E_\theta} f_N (r\theta) d\theta\\
 & \leqslant \int^{2\pi}_0 f_N (r \theta) d\theta \leqslant c
\endalign
$$
Thus,
$$
|\cup I_\sigma| \leqslant c r^{-(m+1)},
$$
and in particular,
$$
|I_\sigma| \leqslant c r^{-(m+1)}.
$$
Thus on the circle $|t| = r$, the arc $rI_\sigma$, has measure
$r^{-m}$.  Now at the end-points of arc $rI_\sigma$ on the circle
$|t|=r, f_N(r\theta) \leqslant r^{m+1}$.  Since $f_N$ satisfies (3.5),
we may apply Harnack's inequality, Lemma (3.3) to $f_N$ on $rI_\sigma$.  The potential in
(3.5) on $rI_\sigma$ is bounded by $r^{2m}$, and $rI_\sigma$ has
length at most $r^{-m}$, so (3.2) applies.  Thus,
$$
\underset{\theta \in I_\sigma}\to{\sup} f(r\theta) \leqslant C r^{m+1}.
$$
Thus we have our conclusion for $r\geqslant 100$.  For $r\leqslant
100$ since Harnack's inequality applies to (3.5), and since
$f_N(0)=1$, and $\lambda_N$ bounded, the conclusion of our lemma is
again immediate.
\enddemo

$\hfill \blacksquare$

For $\Bbb R^n$, we define, for fixed $r, |t|=r$ as before 
$$
E_\omega = \{ \omega \in S^{n-1} : f_N (r\omega) > (1 +
r)^{(m+1)(n-1)}\}.
$$
The proof then follows the $\Bbb R^2$ case, but now using the
properties of Whitney cubes which replace the arcs $I_\sigma$.  That
is we can cover $S^{n-1}$ by coordinate patches, such that $S^{n-1}$
is approximately Euclidean on each patch.  We then decompose
$E_\omega$ into $n-1$ dimensional cubes on $S^{n-1}$ using the Whitney
lemma [St].

\demo{Proof of Theorem 3.2}  The uniform bounds of Lemma (3.8) and the
bounds on $\lambda_N$ from Lemma (3.5) allow us to take a limit of a
suitable subsequence of $N$, with $N\to \infty$.  Since $f_N(0) = 1$,
for all $N$, we get by standard elliptic estimates a $\lambda, 0 <
\lambda < \infty$, and a function $f(t), f(0) = 1, f > 0$, which
satisfies the Schr\"odinger equation of our theorem.
\enddemo

$\hfill\blacksquare$
\medskip
The next proposition handles the case of equality in (3.1) when $s=0$.

\proclaim{Proposition 3.9} Assume that $s=0,\, n=2$ and we have equality
in (3.1). Then for every $\lambda>0$ we can construct an oscillatory
and bounded solution $f(t)$ with $f(0)=1$, or we have a solution
$f(t)$ that satisfies $f(0)=1$ and for some $C>0$
  $$ 0<f(t)<e^{C|t|}.$$
\endproclaim
\demo{Proof} Since $s=0$ and we have equality in (3.1), necessarily
$j=2m$. Thus $q(t)$ vanishes identically or after a rotation of
coordinates $q(t)=t_1^{2m}$. If $q(t)$ vanishes identically
we have the spectral problem
$$ \Delta f +\lambda f=0$$
We select $f(t)=\cos \lambda^{1/2} t_1$ as a solution which has the claimed
properties. In case $q(t)=t_1^{2m}$ we may separate variables. Let $\psi(z)$
denote the first eigenfunction with eigenvalue $\lambda_0$ of the problem
   $$-\psi''(z)+z^{2m}\psi(z)=\lambda_0\psi,\, -\infty<z<\infty.$$
We are seeking a solution to the spectral problem
$$-\Delta f+(t_1^{2m}-\lambda)f=0$$
If $\lambda>\lambda_0$, we take $f(t)=\psi(t_1)\cos \mu t_2$, 
$\mu =|\lambda-\lambda_0|^{1/2}$ as a solution. If 
$\lambda\leqslant \lambda_0$ we take $f(t)=\psi(t_1)e^{\mu t_2}$ 
as a solution. This proves the proposition.
\enddemo
$\hfill\blacksquare$

We may combine the previous proposition and theorem (3.2) to give
a complete solution to the spectral problem when $s=0$ and $n=2$.

\proclaim{Proposition 3.10} Assume that $q(t),\, t\in \Bbb R^2$ is 
a non-negative
homogeneous polynomial of degree $2m$ with real coefficients. Then
for some $\lambda>0$ there exists in all $\Bbb R^2$ a solution $f(t)$
to the spectral problem
$$-\Delta f+(q(t)-\lambda)f=0,$$
where $f(0)=1$ and $f(t)$ satisfies with some $C>0$ the growth condition
  $$ |f(t)|\leqslant e^{C|t|}$$
which is the growth condition which satisfies the requirement (2.20)
for $s=0$.
\endproclaim

We are now in a position to combine all of our results and state a 
concrete theorem
about analytic hypoellipticity

\proclaim{Theorem 3.11} Assume that the basis $\Cal B$ for $V_1$
contains an element from the center of $\frak n$. Assume that $\#\Cal S=2$
and the elements of $\Cal S$ commute amongst themselves. Then the
sums of squares operator (1.1) has a solution $u$ that is Gevrey
of order exactly $m+1$ at the origin.
\endproclaim

\demo{Proof} Theorem (2.5) and proposition (2.2) reduces the spectral
problem to that considered in proposition (3.10). Proposition
(3.10) guarantees for us a solution with the requisite growth
properties that allows us to apply Theorem (2.3). Theorem (2.3)
allows us to conclude the existence of the solution $u$ with the given 
Gevrey property.
\enddemo
$\hfill\blacksquare$

\medskip
We believe Theorem (3.2) holds in full generality without the
assumption (3.1) and with the weaker conclusion (2.20) which is
enough.  The following unpublished result of F. Treves seems to
support this guess.
\proclaim{Theorem 3.12  (F. Treves) }  Let $p(t)$, be any homogeneous
polynomial on $\Bbb R^n$ of degree $2s$, with real coefficients.  Then
there exist a solution $f(t)$ to,
$$
\Delta f+ p (t) f(t) = 0,
$$
in all $\Bbb R^n$, such that, for some $C_1, C_2 > 0$,
$$
|f(t) | \leqslant C_1 \exp ( C_2 |t|^{s+1}).
$$
\endproclaim

\medskip

\subhead \S 4. Examples and further remarks\endsubhead

\noindent{\bf Example (4.1) }  We consider
$$
\align
X_1 &= \frac{\partial}{\partial x_1} , X_2 = \frac{\partial}{\partial
x_2}, X_3 = (x_1^2 + x^2_2)\frac{\partial}{\partial x_3}\\
  & X_4 = x^7_1 x^2_2 \frac{\partial}{\partial x_4}, X_5 = x_1^2 x^7_2
\frac{\partial}{\partial x_5}.
\endalign
$$
Let $L= \sum^5_{i=1} X_i^2$.  We define,
$$
u(x_1, x_2, x_3, x_4, x_5) = \int^\infty_0 f(\rho x_1, \rho x_2)
e^{\lambda^{1/2} \rho^3 x_3 + i \rho^{10} (x_4 + x_5)} e^{-M\rho^3}
dp.\tag 4.1
$$
$u$ will satisfy $Lu=0$, provided $f$ satisfies in $\Bbb R^2$, 
$$
-\Delta f + \left[(x_1^{14} x_2^4 + x_1^4 x_2^{14}) - \lambda (x_1^2 +
x^2_2)^2\right] f=0\tag 4.2
$$
In (4.2), $m=9, s = 2$ and $j=4$.  Thus (3.1) is satisfied.  Thus
Theorem (3.2) applies and we have a solution $f(x_1, x_2), f(0) = 1$,
with $x=(x_1, x_2)$,
$$
0< f (x) \leqslant C(1+|x|)^{10} \tag 4.3
$$
Thus the integral in (4.10) is convergent.

Further,
$$
\frac{\partial^\sigma u(0)}{\partial x^\sigma_4} =
i^\sigma\int^\infty_0 \rho^{10\sigma} e^{-M\rho^3} dp.
$$
Thus the solution is exactly in Gevrey $10/3 = (m+1)/(s+1)$ at the origin.

Note that $\xi_3=0$, defines the non-symplectic stratum,
corresponding to $\frac{\partial}{\partial x_3}$.  Now if one were to
consider the Lie Algebra $\frak n$, with the same bracket relations as
the $X_i$ above, one finds,
$$
\frak n = \overset{10}\to{\underset{i=1}\to\oplus} V_i.
$$  
The co-vector corresponding to $\frac{\partial}{\partial x_3}$ when 
restricted to
$V_3$ does not vanish, but vanishes when restricted to $V_i,\, i\neq 3$ 
i.e. $s=2$.  Next $\dim V_{10} = 2$, corresponding to
$\frac{\partial}{\partial x_4}$ and $\frac{\partial}{\partial x_5}$, and $m=9$
clearly.  Note too the strong form of assumption (2.4) holds, $[V_i,
V_j] = 0, i, j > 1$.  Next Remark (2.4) is also clear.  We have the
right to induce representations using
$$
\lambda_2 = (0, 0, \dots, 0, \lambda_{10,1} , \lambda_{10,2})
$$
where $(0, \dots, 0, \lambda_{10,1}, 0)$ is a co-vector to
$\frac{\partial}{\partial x_4}$ and $(0, \dots, 0, 0, \lambda_{10,2})$ is a
co-vector to $\frac{\partial}{\partial x_5}$.  Given this flexibility,
one may desire to set $\lambda_{10, 1} \neq 0$ but $\lambda_{10, 2} =
0$.  This is equivalent to searching for a solution $u$ of the form
$u(x_1, x_2, x_3, x_4)$ and we are led to find an $f(x_1, x_2)$ a
solution to,
$$
-\Delta f + \left[ x_1^{14} x_2^4 - \lambda (x^2_1 + x_2^2)^2\right] f
=0. \tag 4.4 
$$
This operator is less positive then (4.2) and Theorem (3.2) does not
apply to (4.4).  Thus to apply Theorem (3.2) we must keep both
$\lambda_{10, 1} , \lambda_{10, 2}$ non-zero and set them equal to 1.
Thus in the parameter $\lambda_{m+1, \beta}$ there is great
flexibility, and given that we have not proved Theorem (3.2) in full
strength it is advisable to induce in all the parameters
$\lambda_{m+1,\beta}$ so that Theorem (3.2) applies.  Further $\Cal S
= \{\frac{\partial}{\partial x_1} , \frac{\partial}{\partial x_2}\},
\# \Cal S =2$, and the elements of $\Cal S$ commute.
Notice also $\sum^3_{i=1} X^2_i$ has a characteristic set that is symplectic.
\medskip

\noindent{\bf Lifting: }  We have seen in Example (4.1) that the
results of Section 3 apply to even degenerate vector fields.  However,
we have been unable to adapt the lifting mechanism in [RS] to Section
2.  The difficulty is to verify formula (2.23) for the degenerate
case, though (2.23) and (4.1) are very close.  Furthermore, lifting
may introduce non-symplectic strata where non-existed.

\noindent{\bf Example (4.2) } 
$$
X_1 =\frac{\partial}{\partial x_1} , X_2 = x^2_1
\frac{\partial}{\partial x_2}.
$$
The characteristic set of $X_1, X_2$ is symplectic.  The corresponding
lifted vector fields are,
$$
Y_1 = \frac{\partial}{\partial x_1} , Y_2 = \frac{\partial}{\partial
x_3} + x_1 \frac{\partial}{\partial x_4} + x_1^2
\frac{\partial}{\partial x_2}.
$$
The vanishing of $\xi_4 = 0$, the co-vector corresponding to
$\frac{\partial}{\partial x_4}$, defines a non-symplectic stratum, and
$L=Y^2_1 + Y^2_2$ is in fact non-analytic hypoelliptic.  In fact the
Gelfand-Kirillov dimension is 1, and this operator gives an ODE as in
[CC].  In fact by considering a solution $u$ of the form $u=u(x_1,x_2,x_4)$
to $L$ we easily reduce matters to the ODE considered in [CC]. 
However, it is also clear that if there are
non-symplectic strata, then after lifting there will be non-symplectic
strata, and thus the co-vectors in the degenerate fields determine the
parameters in $\frak n^*$ for the induced representation in an
invariant way.

\enddocument